\documentclass[a4paper, 10pt, conference]{ieeeconf}

\IEEEoverridecommandlockouts

\usepackage{amsfonts}
\usepackage{amsmath}
\usepackage{amssymb}
\usepackage{makeidx}
\usepackage{epsfig}
\usepackage{graphicx}

\def\R{\mathcal{R}}
\def\P{\mathcal{P}}
\def\O{\mathcal{O}}
\def\L{\mathcal{L}}
\def\S{\mathcal{S}}
\def\W{\mathcal{W}}

\def\GO{G_\O}
\def\GR{G_\R}
\newtheorem{theorem}{Theorem}
\newtheorem{corollary}[theorem]{Corollary}
\newtheorem{definition}{Definition}

\newtheorem{lemma}[theorem]{Lemma}
\newtheorem{proposition}[theorem]{Proposition}

\begin{document}

\title{Link Failure Detection in Multi-hop Control Networks}

\author{Alessandro D'Innocenzo, Maria Domenica Di Benedetto and Emmanuele Serra
\thanks{The authors are with the Department of Electrical and Information Engineering, University of L'Aquila. Address: Via G. Gronchi, 18 Nucleo Industriale di Pile, L'Aquila, 67100 Italy. Tel: +39 328 941 5922. Email: \{mariadomenica.dibenedetto, alessandro.dinnocenzo, emmanuele.serra\}@univaq.it. The research leading to these results has received funding from the European Union Seventh Framework Programme [FP7/2007-2013]  under grant agreement n°257462 HYCON2 Network of excellence.}}

\maketitle

\vspace{-1.50cm}

\begin{abstract}
A Multi-hop Control Network (MCN) consists of a plant where the communication between sensors, actuators and computational unit is supported by a wireless multi-hop communication network, and data flow is performed using scheduling and routing of sensing and actuation data. We characterize the problem of detecting the failure of links of the radio connectivity graph and provide necessary and sufficient conditions on the plant dynamics and on the communication protocol. We also provide a methodology to \emph{explicitly} design the network topology, scheduling and routing of a communication protocol in order to satisfy the above conditions.
\end{abstract}



\section{Introduction} \label{secIntro}

Wireless networked control systems are spatially distributed control systems where the communication between sensors, actuators, and computational units is supported by a shared wireless communication network. Control with wireless technologies typically involves multiple communication hops for conveying information from sensors to the controller and from the controller to actuators. The use of wireless networked control systems in industrial automation results in flexible architectures and generally reduces installation, debugging, diagnostic and maintenance costs with respect to wired networks. The main motivation for studying such systems is the emerging use of wireless technologies in control systems (see e.g.,~\cite{akyildiz_wireless_2004}, \cite{song_wirelesshart:_2008}, and \cite{song_complete_2008}).

Although Multi-hop Control Networks (MCNs) offer many advantages, their use for control is a challenge when one has to take into account the joint dynamics of the plant and of the communication protocol. Wide deployment of wireless industrial automation requires substantial progress in wireless transmission, networking and control, in order to provide formal models and verification/design methodologies for wireless networked control system. The design of the control system has to consider the presence of the network, as it represents the interconnection between the plant and the controller, and thus affects the dynamical behavior of the system. The analysis of stability, performance, and reliability of real implementations of wireless networked control systems requires addressing issues such as scheduling and routing using real communication protocols.

Recently, a huge effort has been made in scientific research on Networked Control Systems (NCSs), see \cite{Zhang2001},~\cite{WalshCSM2001},~\cite{SpecialIssueNCS2004},~\cite{Hespanha2007}, and \cite{HeemelsTAC10}, and references therein for a general overview. However, the literature on NCSs usually does not take into account the non--idealities introduced by scheduling and routing communication protocols of Multi-hop Control Networks. In~\cite{andersson_simulation_2005}, a simulative environment of computer nodes and communication networks interacting with the continuous-time dynamics of the real world is presented. To the best of our knowledge, the only formal model of a Multi-hop Control Network has been presented in \cite{AlurRTAS09,AlurTAC11}, where the modeling and stability verification problem has been addressed for a MIMO LTI plant embedded in a MCN, when the controller is already designed. A mathematical framework has been proposed, that allows modeling the MAC layer (communication scheduling) and the Network layer (routing) of the recently developed wireless industrial control protocols, such as WirelessHART (\texttt{www.hartcomm2.org}) and ISA-100 (\texttt{www.isa.org}).


\begin{figure*}[t]
\begin{center}
\includegraphics[width=1.0\textwidth]{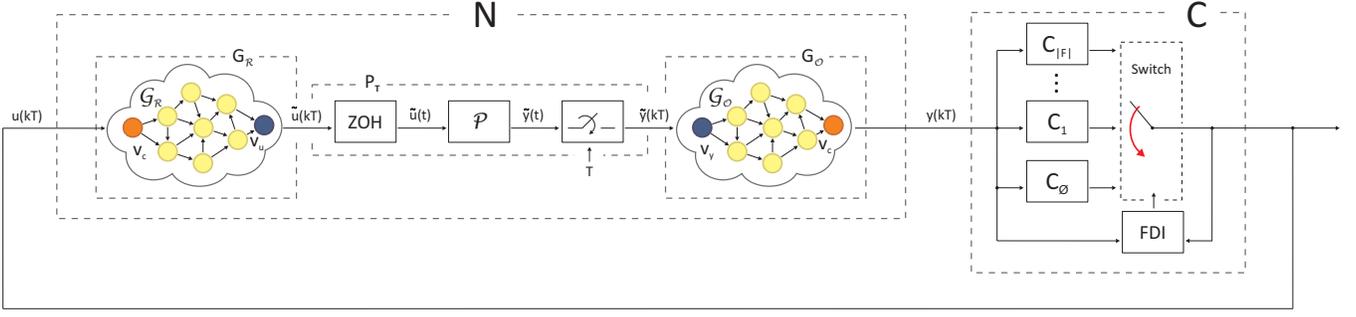}
\caption{Proposed control scheme of a MCN.}\label{fighybridControlScheme}
\end{center}
\end{figure*}

Consider the networked control architecture illustrated in Figure \ref{fighybridControlScheme}, that consists of a plant $\P$ interconnected to a controller $\mathcal C$ via two multi-hop wireless communication networks $\GR$ and $\GO$. We proved in \cite{DiBenedettoIFAC11Stab} that for any \emph{time-invariant} topology $i$ of $\GR$ and $\GO$, characterized by at least one path between the controller and the plant, it is always possible to design a controller $\mathcal C_i$, a routing and a scheduling to arbitrarily assign the eigenvalues of the closed loop system. Consider the following two application scenarios. In the first scenario (e.g. the mine application investigated in \cite{D'InnocenzoCASE2009}), an industrial plant is connected to a controller via a multi-hop wireless communication network: the graph topology of the wireless network is time-varying because of link failures and battery discharge of the communication nodes. In the second scenario, a plant is connected to a controller via a swarm of mobile agents (e.g. robots \cite{ZavlanosCDC07} or UAVs \cite{MeskinTAC2009}) equipped with wireless communication nodes: the graph topology of the wireless network is time-varying because of motion of the agents. In both scenarios, the time-varying topology perturbs the dynamics of the interconnected system $N$, and the controller is required to detect the current topology $i$ of $\GR$ and $\GO$ to apply the corresponding control law $\mathcal C_i$.

In this paper we suppose that the topology of $\GR$ and $\GO$ is \emph{time-varying} because of link failures, and provide a methodology to detect the set of faulty links using Fault Detection and Identification (FDI) methods. In the taxonomy of fault diagnosis techniques, we leverage on the model-based approach introduced by the pioneering works in \cite{BeardPhDThesis1971,JonesPhDThesis1973} on observer-based FDI, later pursued in \cite{MassoumniaTAC89} for linear systems and in~\cite{DePersisTAC2001} for non-linear systems.

As can be inferred from the recent survey \cite{Gupta10}, fault tolerant control and fault diagnosis is one of the main issues addressed in the research on NCSs. However, most of the existing literature on NCSs fault diagnosis (e.g. \cite{WangTSP2008},~\cite{MeskinTAC2009}) usually addresses communication delays, and does not consider the effect of the communication protocol introduced by a Multi-hop Control Network. In \cite{CommaultTAC2007}, a procedure to minimize the number and cost of additional sensors, required to solve the FDI problem for \emph{structured systems}, is presented. In \cite{PappasCDC2010a}, the design of an intrusion detection system is presented for a MCN, where the network \emph{itself} acts as the controller. Our modeling framework differs from that developed in \cite{PappasCDC2010a}, since we model the MCN as an input-output system where the wireless networks \emph{transfer} sensing and actuation data between a plant and a controller (they are \emph{relay} networks), while in \cite{PappasCDC2010a} the MCN is an autonomous system where the wireless network \emph{itself} acts as a controller. Moreover, in our model we explicitly take into account the effect of the scheduling ordering of the node transmissions in the sensing and actuation data relay.

Our work differs from the existing literature since we characterize the communication link failures detection problem in a MCN as a FDI problem, and state necessary and sufficient conditions \emph{on the plant dynamics} and \emph{on the communication protocol}. Moreover, we provide a methodology to \emph{explicitly} design the network topology, scheduling and routing of a communication protocol in order to satisfy link failure detection conditions of a MCN for any failure of communication links. The explicit design of scheduling and routing is a fundamental aspect of our contribution. In fact, as evidenced in \cite{D'InnocenzoCASE2009}, when applying a wireless industrial control protocol to the real scenario the topology of the wireless network introduces hard limitations in the choice of the scheduling. This is due to the fact that most of the wireless industrial control protocols suggest that the communication scheduling satisfies a specific ordering (see \cite{D'InnocenzoCASE2009,DiBenedettoIFAC11Rout} for more details). The results in \cite{DiBenedettoIFAC11Stab} and in this paper mitigate these constraints, by proving that it is not required to perform scheduling according to a specific ordering. This allows to strongly reduce the scheduling length, as illustrated in \cite{DiBenedettoIFAC11Stab}.


\section{Modeling of MCNs} \label{secMCNModeling}

\begin{figure}[ht]
\begin{center}
\includegraphics[width=0.5\textwidth]{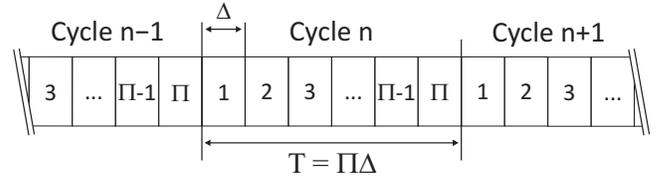}
\caption{Time-slotted structure of frames.}\label{frame}
\end{center}
\end{figure}

The challenges in modeling MCNs are best explained by considering the recently developed wireless industrial control protocols, such as WirelessHART and ISA-100. These standards require that designers of wireless control networks define a communication scheduling for all communication nodes of a wireless network. For each working frequency, time is divided into slots of fixed duration $\Delta$, and groups of $\Pi$ time slots are called frames of duration $T = \Pi \Delta$ (see Figure \ref{frame}). For each frame, a communication scheduling allows each node to transmit data only in a specified time slot and frequency, i.e. a mixed TDMA and FDMA MAC protocol is used. The communication scheduling is periodic with period $\Pi$, i.e. it is repeated in all frames.
The standard specifies a syntax for defining scheduling and routing and a mechanism to apply them, but the issue of designing them remains a challenge for engineers and is currently done using heuristic rules. To allow systematic methods for designing the communication protocol configuration, a mathematical model of the effect of scheduling and routing on the control system is needed.
\begin{definition}\label{defMCN}
A SISO Multi-hop Control Network is a tuple $N = ( \P, G_\R, \eta_{\R}, G_\O, \eta_{\O}, \Delta)$ where:
\begin{itemize}
\item $\P = (A^c_{\P},B^c_{\P},C^c_{\P})$ models a plant dynamics in terms of matrices of a continuous-time SISO LTI system.

\item $G_\R = ( V_{\R},E_{\R}, W_{\R} )$ is the controllability radio connectivity acyclic graph, where the vertices correspond to the nodes of the network, and an edge from $v$ to $v'$ means that $v'$ can receive messages transmitted by $v$ through the wireless communication link $(v,v')$. We denote $v_{c}$ the special node of $V_{\R}$ that corresponds to the controller, and $v_u \in V_{\R}$ the special node that corresponds to the actuator of the input $u$ of $\P$. The weight function $W_{\R} : E_{\R} \to \mathbb R^+$ associates to each link a positive constant. The role of $W_{\R}$ will be clear in the following definition of $\eta_{\R}$.

\item $\eta_{\R} \colon \mathbb N \to 2^{E_{\R}}$ is the controllability communication scheduling function, that associates to each time slot of each frame a set of edges of the controllability radio connectivity graph. Since in this paper we only consider a periodic scheduling that is repeated in all frames, we define the controllability communication scheduling function by $\eta_{\R} \colon \{1, \ldots, \Pi\} \to 2^{E_{\R}}$. The integer constant $\Pi$ is the period of the controllability communication scheduling. The semantics of $\eta_{\R}$ is that $( v,v' ) \in \eta(h)$ if and only if at time slot $h$ of each frame the data content of the node $v$ is transmitted to the node $v'$, multiplied by the weight $W_{\R}(v,v')$. We assume that each link can be scheduled only one time for each frame. This does not lead to loss of generality, since it is always possible to obtain an equivalent model that satisfies this constraint by appropriately splitting the nodes of the graph, as already illustrated in the memory slot graph definition of \cite{AlurTAC11}.

\item $G_\O = ( V_{\O},E_{\O}, W_{\O} )$ is the observability radio connectivity acyclic graph, and is defined similarly to $G_\R$. We denote with $v_{c}$ the special node of $V_{\O}$ that corresponds to the controller, and $v_y \in V_{\O}$ the special node that corresponds to the sensor of the output $y$ of $\P$.

\item $\eta_{\O} \colon \{1, \ldots, \Pi\} \to 2^{E_{\O}}$ is the observability communication scheduling function, and is defined similarly to $\eta_{\R}$. We remark that $\Pi$ is the same period as the controllability scheduling period.

\item $\Delta$ is the time slot duration. As a consequence, ${T} = \Pi \Delta$ is the frame duration.
\end{itemize}
\end{definition}
Definition \ref{defMCN} allows modeling communication protocols that specify TDMA, FDMA and/or CDMA access to a shared communication resource, for a set of communication nodes interconnected by an arbitrary radio connectivity graph. In particular, it allows modeling wireless multi-hop communication networks that implement protocols such as WirelessHART and ISA-100. Our MCN model differs from the framework developed in \cite{AlurTAC11}, since it allows modeling redundancy in data communication sending control data through multiple paths in the same frame and then merging these components according to the weight function. This kind of redundancy is called \emph{multi-path routing} (or \emph{flooding}, in the \emph{communication} scientific community), and aims at rendering the MCN robust with respect to link failures and to mitigating the effect of packet losses.

For any given radio connectivity graph that models the communication range of each node, designing a scheduling function induces a communication scheduling (namely the time slot when each node is allowed to transmit) and a multi-path routing (namely the set of paths that convey data from the input to the output of the connectivity graph) of the communication protocol. Since the scheduling function is periodic the induced communication scheduling is periodic, and the induced multi-path routing is static.

We define a connectivity property of the controllability and observability graphs with respect to the corresponding scheduling.
\begin{definition}
Given a controllability graph $G_{\R}$ and scheduling $\eta_{\R}$, we define $G_{\R}(\eta_{\R}(h))$ the sub-graph of $G_{\R}$ induced by keeping the edges scheduled in the time slot $h$. We define $G_{\R}(\eta_{\R}) = \bigcup\limits_{h = 1}^{\Pi} G_{\R}(\eta_{\R}(h))$ the sub-graph of $G_{\R}$ induced by keeping the union of edges scheduled during the whole frame.
\end{definition}
\begin{definition}
We say that a controllability graph $G_{\R}$ is jointly connected by a controllability scheduling $\eta_{\R}$ if and only if there exists a path from the controller node $v_c$ to the actuator node $v_u$ in $G_{\R}(\eta_{\R})$.
\end{definition}
The above definitions can be given similarly for observability graph $G_{\O}$ and scheduling $\eta_{\O}$.

The dynamics of a MCN $N$ can be modeled by the interconnection of blocks as in Figure \ref{fighybridControlScheme}. The block $P_{T}$ is characterized by the discrete-time state space representation $(A_{\P},B_{\P},C_{\P})$ obtained by discretizing $(A^c_{\P},B^c_{\P},C^c_{\P})$ with sampling time ${T} = \Pi \Delta$. We assume that the plant $\P$ is stabilizable and detectable, and that $\P = (A^c_{\P},B^c_{\P},C^c_{\P})$ is the controllable and observable minimal representation. If this assumption does not hold, then even with an ideal interconnection between the controller and the plant it is clearly not possible to stabilize the closed loop system, and the control scheme in Figure \ref{fighybridControlScheme} looses any interest.

The block $G_{\R}$ models the dynamics introduced by the data flow of the actuation data through the communication network represented by $G_\R$ according to the applied controllability scheduling $\eta_{\R}$. In order to define the dynamical behavior of $G_{\R}$, we need to define the dynamics of the data flow through the network, according to the scheduling $\eta_{\R}$.

We associate to the controller node $v_c$ a real value $\mu_c(kT)$ at time $k$, and we assume that $v_c$ is periodically updated with a new control command at the beginning of each frame and holds this value for the whole duration of the frame. Formally, $\mu_c(kT) = u(kT)$.

The dynamics of the other nodes needs to be defined at the level of time slots. We associate to each other node $v_j~\in~V_{\R} \setminus \{v_c\}$ a real value $\mu_{i,j}(h)$ at time slot $h$ for each node $v_i$ belonging to the set $inc(v_j) = \{v \in V_{\R} : (v, v_j) \in E_{\R}\}$ of edges incoming in $v_j$.

When the link from $v_i$ to $v_j$ is not scheduled at time slot $h$, the variable $\mu_{i,j}(h)$ is not updated. When the link from $v_i$ to $v_j$ is scheduled at time slot $h$, the variable $\mu_{i,j}(h)$ is updated with the sum of the variables associated to node $v_i$ in the time slot $h$ multiplied by the link weight $W_{\R}(v_i,v_j)$. Formally, for each $v_j \in V_{\R} \setminus \{v_c\}$ and for each time slot $h \in \{1, \ldots, \Pi\}$:
\begin{equation*}
\mu_{i,j}(h+1) = \left\{
\begin{array}{l}
  \mu_{i,j}(h) \text{ if } (v_i, v_j) \notin \eta_{\R}(h), \\ \\
  W_{\R}(v_i,v_j) \cdot \sum_{v_k \in inc(v_i)} \mu_{k,i}(h) \\
  \text{if } (v_i, v_j) \in \eta_{\R}(h).
\end{array}
\right.
\end{equation*}
Finally, the actuator node $v_u$ periodically actuates a new actuation command at the beginning of each frame on the basis of its variables $\mu_{i,u}$, and holds this value for the whole duration of the frame. Formally,
\begin{equation*}
\tilde u(kT) = \sum_{v_i \in inc(v_u)} \mu_{i,u}(kT).
\end{equation*}
The following proposition proved in \cite{DiBenedettoIFAC11Stab} characterizes the dynamics of $G_{\R}$ at the level of frames, induced by the data flow through the network at the level of time slots.
\begin{proposition}\label{propGz}
\cite{DiBenedettoIFAC11Stab} Given $G_{\R}$ and $\eta_{\R}$, the controllability graph can be modeled as a discrete time SISO LTI system with sampling time equal to the frame duration $T = \Pi \Delta$, and characterized by the following transfer function:
\begin{equation*}
G_{\R}(z) = \sum_{d = 1}^{D_{\R}} \frac{\gamma_{\R}(d)}{z^d},
\end{equation*}
where $D_{\R} \in \mathbb N$ is the maximum delay introduced by $G_{R}$, and $\forall d \in \{1, \ldots, D_\R -1\}$, $\gamma_{\R}(d) \in \mathbb R_0^+$, $\gamma_{\R}(D_{\R}) \neq 0$.
\end{proposition}

\begin{figure}[ht]
\begin{center}
\includegraphics[width=0.45\textwidth]{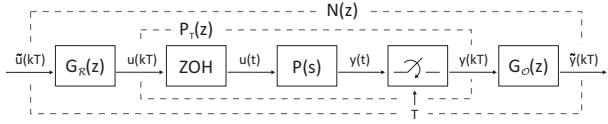}
\caption{Transfer function of the MCN interconnected system.} \label{MCNblocks}
\end{center}
\end{figure}
$G_{\O}(z)$ can be computed similarly. The dynamics of a MCN $N$ can be modeled as in Figure \ref{MCNblocks}, where each block is a discrete time SISO LTI system with sampling time equal to the frame duration, characterized by the transfer functions $G_{\R}(z)$, $P_T(z)$ and $G_{\O}(z)$.

Let $x_{\O} \in \mathbb R^{n_{\O}}$, $x_{\P} \in \mathbb R^{n_{\P}}$ and $x_{\R} \in \mathbb R^{n_{\R}}$ be respectively the states of the observability graph, of the plant, and of the controllability graph. We will denote by $x = \left[\begin{array}{ccc}
                                                                                           x_{\O}^\top & x_{\P}^\top & x_{\R}^\top
                                                                                         \end{array}\right]^\top$
the extended state of $N$, with $x \in \mathbb R^{n}$, and $n = n_{\O} + n_{\P} + n_{\R}$. The dynamics of $N$ can also be described by the following state space representation:
\begin{align}\label{eqNominalModel}
&x((k+1)T) = A x(kT) + B {u}(kT), \qquad {y}(kT) = C x(kT),\notag \\
&u(kT), \ y(kT) \in \mathbb{R},
\end{align}
with:

\begin{align*}
&A = \left[
                    \begin{array}{lll}
                      A_{\O} & B_{\O}C_{\P} & \textbf{0}_{n_\O \times n_\R} \\
                      \textbf{0}_{n_\P \times n_\O} & {A_{\P}} & B_{\P} C_{\R} \\
                      \textbf{0}_{n_\R \times n_\O} & \textbf{0}_{n_\R \times n_\P} & A_{\R} \\
                    \end{array}
                  \right]\mbox{, }\\
&B = \left[
    \begin{array}{l}
        \textbf{0}_{n_\O \times 1} \\
        \textbf{0}_{n_\P \times 1} \\
        B_{\R} \\
    \end{array}
\right]\mbox{, }
C = \left[
    \begin{array}{l}
        C_{\O}^\top\\
        \textbf{0}_{n_\P \times 1} \\
        \textbf{0}_{n_\R \times 1} \\
    \end{array}
\right]^\top,
\end{align*}
and

{\small

\begin{align*}
&A_\R = \left[
          \begin{array}{lll}
            0 &\ \gamma_{\R}(D_\R) &\ \ \gamma_{\R}(D_\R - 1) \cdots \gamma_{\R}(2)\\
            \textbf{0}_{(D_\R-2) \times 1} &\ \textbf{0}_{(D_\R-2) \times 1} &\ \ \textbf{I}_{D_\R -2}\\
            0 &\ 0 &\ \ \textbf{0}_{1 \times (D_\R-2)}
          \end{array}
        \right],\\
&B_\R = \left[
          \begin{array}{l}
            \gamma_{\R}(1) \\
            \textbf{0}_{1 \times (D_\R-2)} \\
            1
          \end{array}
        \right],
C_\R = \left[
          \begin{array}{l}
            1 \\
            \textbf{0}_{(D_\R -1) \times 1}
          \end{array}
        \right]^\top. 
\end{align*}

}

The matrices $(A_{\O},B_{\O},C_{\O})$ are defined similarly.

\section{Fault Detection on MCNs} \label{secFailureModelMCN}

In this section we provide a methodology to detect the current dynamics of a MCN subject to link failures using Fault Detection and Identification (FDI) methods. The failure of a set of links $f \subseteq {E_{\R} \cup E_{\O}}$ on the dynamics \eqref{eqNominalModel} can be modeled as follows:
\begin{eqnarray}\label{eqFaultyMCNDynamics}
x((k+1)T) &=& A x(kT) + B u(kT) + L_{f} m_{f}(kT) \notag \\
y(kT) &=& C x(kT)
\end{eqnarray}
where $m_{f}(kT) : \mathbb N \to \mathbb{R}^{n+1}$ is an arbitrary function of time and $L_{f} \colon \mathbb{R}^{n+1} \to \mathbb{R}^{n}$ is called the failure signature map associated to the configuration of failures $f$. We define the failure signature maps as in Figure \ref{figL}:

\begin{figure*}[t]
\begin{center}
\includegraphics[width=10px]{blank}
\end{center}
\end{figure*}

\begin{figure*}[t]
\begin{picture}(0,0)
\put(70,10){
\scalebox{1}{
$
L_{f} = \left[
          \begin{array}{llll}
            0 & -\delta_{\O, f} & \textbf{0}_{1 \times n_\P} &  \textbf{0}_{1 \times n_\R} \\
            \textbf{0}_{(n_\O + n_\P -1) \times 1} & \textbf{0}_{(n_\O + n_\P -1) \times n_\O} & \textbf{0}_{(n_\O + n_\P -1) \times n_\P} &  \textbf{0}_{(n_\O + n_\P -1) \times n_\R} \\
            0 & \textbf{0}_{1 \times n_\O} & \textbf{0}_{1 \times n_\P} &  -\delta_{\R, f}\\
            \textbf{0}_{(n_\R -1) \times 1} & \textbf{0}_{(n_\R -1) \times n_\O} & \textbf{0}_{(n_\R -1) \times n_\P} & \textbf{0}_{(n_\R -1) \times n_\R} \\
          \end{array}
        \right],
$
}
}
\end{picture}
\vspace{0.5cm}
\caption{Matrix $L_{f}$.}
\label{figL}
\end{figure*}

where the $d$-th components $\delta_{\R,f}(d)$ and $\delta_{\O,f}(d)$ of the row vectors $\delta_{\R, f}~=~\left[\begin{array}{ccc}
                                              \delta_{\R, f}(D_\R) & \cdots & \delta_{\R, f}(1)
                                            \end{array}\right]$
and $\delta_{\O, f}~=~\left[\begin{array}{ccc}
                                              \delta_{\O, f}(D_\O) & \cdots & \delta_{\O, f}(1)
                                            \end{array}\right]$
are the perturbations introduced by the configuration of failures $f$ in the paths of $G_\R$ and $G_\O$ characterized by delay $d$. Since $\gamma_{\R}(d) \geq 0$ and $\gamma_{\O}(d) \geq 0$, and a failure of each path reduces the value of the corresponding component, then $\delta_{\R, f}(d) \geq 0$ and $\delta_{\O, f}(d) \geq 0$ for each $f \subseteq {E_{\R} \cup E_{\O}}$. In the absence of failures $L_{\varnothing} = \textbf{0}_{n \times (n+1)}$.

The signal $m_{f}(kT)$ depends on the protocol applied by the communication nodes when the configuration of failures $f$ occurs. By an appropriate choice of $m_{f}(kT)$, it is possible to model by \eqref{eqFaultyMCNDynamics} the dynamics of $N$ when a failure occurs in the set of links $f$, for any protocol applied by the communication nodes in case of failure. As an example, if a node sets to 0 the data contribution incoming from a faulty link, then we can model this behavior by defining $m_{f}(kT)~=~\left[
\begin{array}{cc}
x(kT)^\top & u(kT)^\top
\end{array}\right]^\top$. If a node uses the latest data received from a faulty link, then we can model this behavior by defining $m_{f}(kT)~=~\left[
\begin{array}{cc}
x(kT)^\top & u(kT)^\top
\end{array}\right]^\top + \nu$, with $\nu \in {\mathbb R}^{n+1}$ a constant vector of real numbers.

To perform failure detection of a MCN with the aim of applying an appropriate control law for each dynamics induced by all failure configurations, we first need to define the set $\Phi \subseteq 2^{E_\R \cup E_\O}$ of failures we are interested in distinguishing. In fact, we need to distinguish two failures induced by sets of links $f$, $f'$ only when they introduce different perturbations of the dynamics \eqref{eqNominalModel}, namely when $L_{f} m_f(kT) \neq L_{f'} m_{f'}(kT)$. For this reason, we define $\Phi_\Omega$ the set of equivalence classes $[f]$, each consisting of sets of links that affect the dynamics \eqref{eqNominalModel} by means of the same representative failure signal $L_f m_f(kT)$:
$$
[f] = \left\{f' \subseteq {E_\R \cup E_\O} : \forall k \geq 0, L_{f'} m_{f'}(kT) = L_f m_f(kT) \right\}.
$$
For simplicity of notation, we will denote in the following the equivalence class $[f]$ by a representative set of links $\varphi \in [f]$.
In order to take into account simultaneous failures, we define the subset $\Phi_\Sigma \subset \Phi_\Omega$ of equivalence classes such that the perturbation introduced can be obtained as the sum of perturbations introduced by equivalence classes of $\Phi_\Omega$:
\begin{align*}
\Phi_{\Sigma} = \Bigg\{&f \in \Phi_\Omega : \Big(\exists\ p \in \mathbb N, \exists\ f_1, \ldots, f_p \in \Phi_\Omega \setminus f :\\
&L_f m_f(kT) = \sum_{i=1}^m  L_{f_i} m_{f_i}(kT)\Big)\Bigg\}.
\end{align*}

Define the set of failures as $\Phi = \Phi_\Omega \setminus \Phi_{\Sigma}$. $\Phi$ always contains the equivalence class $\varnothing$, that corresponds to the absence of failures. It is easy to prove that the set $\Phi$ always exists and is unique. For this reason, we can associate to any given MCN $N$ the corresponding unique set of failures $\Phi$ we are interested in distinguishing, and model their simultaneous occurrence as follows:
\begin{align}\label{eqFaultyMCNDynamicsSimultaneousFailures}
&x((k+1)T) = A x(kT) + B u(kT) + \sum\limits_{\varphi \in \Phi} L_{\varphi} m_{\varphi}(kT),\notag\\
&y(kT) = C x(kT).
\end{align}

Given a MCN $N$ and the corresponding faulty set $\Phi$ modeled by~\eqref{eqFaultyMCNDynamicsSimultaneousFailures}, we address the problem of detecting a failure $\varphi \in \Phi$ that is perturbing the dynamics of $N$ by using the measures of the signals $u(\cdot)$, $y(\cdot)$. To this aim we leverage on the model-based approach developed in \cite{MassoumniaTAC89}, which exploits a bank of LTI observer-like systems (called the residual generators) that take as input the signals $u(\cdot)$, $y(\cdot)$, and provides asymptotic estimates of $m_{\varphi}(kT)$ for any failure $\varphi \in \Phi$. This allows to identify which failures are affecting the dynamics of $N$. The problem of designing such residual generators with arbitrary asymptotic convergence rate on the model \eqref{eqFaultyMCNDynamicsSimultaneousFailures} is well known as the \emph{Extended Fundamental Problem in Residual Generation} (EFPRG). Necessary and sufficient conditions for solving the EFPRG have been stated in \cite{MassoumniaTAC89}:


\begin{theorem}
Given the failure model \eqref{eqFaultyMCNDynamicsSimultaneousFailures}, the EFPRG has a solution for the failure $\varphi \in \Phi$ if and only if:
\begin{equation}\label{conditionExtFPRG}
\S^*(\bar{\L}_{\varphi}) \cap \L_\varphi = \textbf{\emph{0}}, \hspace{10px}
\end{equation}
where $\bar{\L}_{\varphi} := \sum_{\varphi' \in \Phi \setminus \varphi}\L_{\varphi'}$.
\end{theorem}
\smallskip
Given any $\L \subseteq \mathbb R^n$, the computation of $\S^*(\L)$ can be performed by applying the (C,A)-Invariant Subspace Algorithm (CAISA) and the UnObservability Subspace Algorithm (UOSA), recursive algorithms provided in \cite{WonhamLMC1979}. We define $\W^*(\L)$ the fixed point of the following recursion (CAISA):
\begin{equation*}
\W_{k+1}(\L) = \L + A\big(\W_k(\L) \cap \mathcal N(C)\big), \quad \W_0(\L) = \textbf{\emph{0}}.
\end{equation*}
We define $\S^*(\L)$ the fixed point of the following recursion (UOSA):
\begin{equation*}
\S_{k+1}(\L) = \W^*(\L) + A^{-1}\big(\S_k(\L)\big) \cap \mathcal N(C), \quad \S_0(\L) = \mathbb R^n.
\end{equation*}
The following lemma provides a useful property of the CAISA and UOSA Algorithms.
\begin{lemma}\label{lemPropUOSA}
Let $\L \subseteq \mathcal N^{\perp}(C)$, then $\W^*(\L) = \L$, and $\S^*(\L) = \L + \mathcal K$ with $\mathcal K \subseteq \mathcal N(C)$. Moreover, if $\L~=~\big(\mathcal N(C)\big)^\perp$, then $\S^*(\L) = \mathbb R^n$. 

\begin{proof}
Let $\L \subseteq \big(\mathcal N(C)\big)^\perp$, then
\begin{align*}
\W_1(\L) &= \L + A\big(\textbf{\emph{0}} \cap \mathcal N(C)\big) = \L + A(\textbf{\emph{0}}) = \L,\\
\W_2(\L) &= \L + A\big(\L \cap \mathcal N(C)\big) = \L + A(\textbf{\emph{0}}) = \L = \W^*(\L).
\end{align*}
For each $k>0$,
\begin{align*}
\S_{k+1}(\L) = \L + A^{-1}\big(\S_k(\L)\big) \cap \mathcal N(C) = \L + \mathcal K_k,
\end{align*}
with $\mathcal K_k \subseteq \mathcal N(C)$. Moreover, if $\L = \big(\mathcal N(C)\big)^\perp$, then:
\begin{align*}
\S_{1}(\L) &= \L + A^{-1}(\mathbb R^n) \cap \mathcal N(C) = \L + \mathbb R^n \cap \mathcal N(C)\\
&=\L + \mathcal N(C) = \big(\mathcal N(C)\big)^\perp + \mathcal N(C) = \mathbb R^n = \S^{*}(\L).
\end{align*}
\end{proof}
\end{lemma}
\smallskip

For the sake of clarity, we address the link failure detection problem starting by two special cases. In the first case, we consider a multi-hop interconnection between the controller and the actuator and a single-hop interconnection between the sensor and the controller, namely the controllability graph $G_\O$ consists of two nodes connected by one link. In the second case, we consider a single-hop interconnection between the controller and the actuator, namely the controllability graph $G_\R$ consists of two nodes connected by one link, and a multi-hop interconnection between the sensor and the controller. In the third case, we consider the general case when both $G_\R$ and $G_\O$ are multi-hop communication networks.


\subsection{$G_\R$ multi-hop and $G_\O$ single-hop} \label{subsecGRmultihop}

If $G_\O$ consists of a single-hop, then $n_\O~=~1$, $A_{\O} = 0$, $B_{\O} = C_{\O} = 1$. As illustrated in \cite{MassoumniaTAC89}, each $L_\varphi$ can be assumed monic with no loss of generality, since when failures are not present the corresponding components of $m_\varphi(kT)$ are identically zero. For this reason, by an appropriate choice of $m_\varphi(kT)$, we define the $L_\varphi$ in~\eqref{eqFaultyMCNDynamicsSimultaneousFailures} as follows:
\begin{equation*}
L_{\varphi} = \left[
          \begin{array}{l}
            \textbf{0}_{(n_\O + n_\P) \times n_\R} \\
            -\delta_{\varphi}\\
            \textbf{0}_{(n_\R -1) \times n_\R} \\
          \end{array}
        \right],
\end{equation*}
where $\delta_{\varphi} \in (\mathbb{R}_0^+)^{n_\R}$ is a row vector and $L_{\varphi} \colon \mathbb{R}^{n_\R} \to \mathbb{R}^{n}$. The following theorem states a negative result.
\begin{theorem}
Let a MCN $N$ and the corresponding faulty set $\Phi$ be given, where $G_\R$ is multi-hop and $G_\O$ is single-hop. Then the EFPRG can be solved for each $\varphi \in \Phi$ if and only if $|\Phi| \leq 2$.

\begin{proof}
(\textbf{sufficiency}) If $|\Phi|=1$ then $\Phi = \{\varnothing\}$, and failures are not defined. If $|\Phi|=2$ then $\Phi = \{\varnothing, \varphi\}$. Therefore, $\bar \L_{\varphi} = \L_{\varnothing}$ and $\bar \L_{\varnothing} = \L_{\varphi}$. Since $\L_{\varnothing} = \textbf{\emph{0}}$, it is easy to derive that $\S^*(\L_{\varphi}) \cap \L_{\varnothing} = \textbf{\emph{0}}$ and that $\S^*(\L_{\varnothing}) \cap \L_{\varphi} = \textbf{\emph{0}}$.

(\textbf{necessity}) Assume that $|\Phi| > 2$. Note that all the elements of the matrix $L_\varphi$ are zeros, except the $(n_\O + n_\P +1)$-th row. For this reason:
\begin{equation*}
\forall \ \varphi \in \Phi, \mbox{ } \L_{\varphi} = span[\textbf{e}_{n_\O + n_\P +1}] := \L_\R.
\end{equation*}
Thus, for each $\varphi \in \Phi$, $\bar{\L}_{\varphi} = \L_\R$. Since $\bar{\L}_{\varphi} \subseteq \S^*(\bar{\L}_{\varphi})$, for each $\varphi \in \Phi$ the following holds:
\begin{align*}
\S^*\left(\bar{\L}_{\varphi}\right) \cap \L_\varphi = \S^*(\L_\R) \cap \L_\R = \L_\R \neq \textbf{\emph{0}}.
\end{align*}
\end{proof}
\end{theorem}

The above theorem states that if the controllability graph is multi-hop and the observability graph is single-hop, then it is not possible to distinguish failures in a set $\Phi$, unless $\Phi$ is trivial. In the following section, we will show that more can be done if the controllability graph is single-hop and the observability graph is multi-hop.


\subsection{$G_\R$ single-hop and $G_\O$ multi-hop} \label{subsecGOmultihop}

If $G_\R$ consists of a single-hop, then $n_\R = 1$, $A_{\R} = 0$, $B_{\R} = C_{\R} = 1$. Using the same reasoning as in the above section, we can define a set $\Phi$ of equivalence classes of link failures that equally perturb the dynamics \eqref{eqFaultyMCNDynamicsSimultaneousFailures}. Since in this case the failures occur in the observability graph, by an appropriate choice of $m_\varphi(kT)$ we define $L_{\varphi} \colon \mathbb{R}^{n_\O} \to \mathbb{R}^{n}$ the failure signature map associated to the equivalence classes $\varphi \in \Phi$:
\begin{eqnarray}\label{eqLOi}
L_{\varphi} &=& \left[
          \begin{array}{lll}
            -\delta_{\varphi}\\
            \textbf{0}_{(n-1) \times n_\O}\\
          \end{array}
        \right],
\end{eqnarray}
where $\delta_{\varphi} \in (\mathbb{R}_0^+)^{n_\O} $ is a row vector and each component $\delta_{\varphi}(d)$ is the perturbation introduced by a failure ${\varphi}$ in the paths of $G_\O$ characterized by delay $d$. The following theorem motivates an extension of the model \eqref{eqFaultyMCNDynamicsSimultaneousFailures}.
\begin{theorem}
Let a MCN $N$ and the corresponding faulty set $\Phi$ be given, where $G_\R$ is single-hop and $G_\O$ is multi-hop. Then the EFPRG can be solved for each $\varphi \in \Phi$ only if the following condition holds:
\begin{equation*}
d\Big(\big(\mathcal N(C)\big)^\perp\Big) \geq \sum\limits_{\varphi \in \Phi} d(\L_\varphi) := n_{\Phi}.
\end{equation*}
\begin{proof}
Equation \eqref{eqLOi} implies that $\L_\varphi \subseteq \big(\mathcal N(C)\big)^\perp$ for each $\varphi \in \Phi$. Therefore $\sum\limits_{\varphi \in \Phi} \L_{\varphi} \subseteq \big(\mathcal N(C)\big)^\perp$, which implies that:
\begin{equation}\label{conditionSumLinXc}
d\left(\sum\limits_{\varphi \in \Phi} \L_{\varphi} \right) \leq d\Big(\big(\mathcal N(C)\big)^\perp\Big).
\end{equation}
Condition \eqref{conditionExtFPRG} implies that $\forall \ \varphi,\varphi' \in \Phi$, $\L_{\varphi} \cap \L_{\varphi'} = \textbf{\emph{0}}$. Therefore:
\begin{equation}\label{conditionSumDim}
d\left(\sum\limits_{\varphi \in \Phi} \L_{\varphi} \right) = \sum\limits_{\varphi \in \Phi} d\left(\L_{\varphi}\right).
\end{equation}
Applying \eqref{conditionSumDim} to \eqref{conditionSumLinXc} completes the proof.
\end{proof}
\end{theorem}

The above theorem shows that it is not possible to design a residual generator for each $\varphi \in \Phi$ if the rank of the matrix $C$ is smaller than $n_{\Phi}$.
In particular, in system \eqref{eqNominalModel} the rank of $C$ is 1, and $n_{\Phi}$ is equal to 1 only if the set $\Phi$ is trivial, namely it contains the equivalence class $\varnothing$ and just one equivalence class $\varphi$. For this reason, we need to consider a more general model for the observability graph. More precisely, we consider observability graphs characterized by $n_{S}$ terminating nodes $v_{1}, \ldots, v_{n_{S}}$, with $n_{S} \geq n_{\Phi}$. This can be modeled without loss of generality by redefining matrices $A_\O$, $B_\O$ and $C_\O$ as in Figure \ref{figAOBOCO}:

\begin{figure*}[t]
\begin{center}
\includegraphics[width=10px]{blank}
\end{center}
\end{figure*}

\begin{figure*}[t]
\begin{picture}(0,0)
\put(0,10){
\scalebox{0.75}{
$
A_\O=\left[
          \begin{array}{cc|ccc}
            \textbf{0}_{1 \times n_{S}} & \gamma_{1}(D_\O) & \gamma_{1}(D_\O - 1) & \cdots & \gamma_{1}(2)\\
            \vdots & \vdots & \vdots & \ddots & \vdots\\
            \textbf{0}_{1 \times n_{S}} & \gamma_{n_{S}}(D_\O) & \gamma_{n_{S}}(D_\O - 1) & \cdots & \gamma_{n_{S}}(2)\\
            \hline
            \textbf{0}_{(D_\O -2) \times n_{S}} & \textbf{0}_{(D_\O -2) \times 1} & & \textbf{I}_{D_\O -2} & \\
            \textbf{0}_{1 \times n_{S}} & 0 & & \textbf{0}_{1 \times (D_\O -2)} &
          \end{array}
        \right],
B_\O=\left[
          \begin{array}{ccccc}
            \gamma_{1}(1) & \cdots & \gamma_{n_{S}}(1) & \textbf{0}_{1 \times (D_\O -2)} & 1 \\
          \end{array}
        \right]^{\top},
C_\O=\left[
          \begin{array}{cc}
            \textbf{I}_{n_{S}} & \textbf{0}_{n_{S} \times (D_\O - 1)} \\
          \end{array}
        \right].
$
}
}
\end{picture}
\vspace{0.5cm}
\caption{Matrices $A_\O$, $B_\O$ and $C_\O$.}
\label{figAOBOCO}
\end{figure*}

where $n_\O~=~D_\O + n_{S}-1$ is the new dimension of the state space. The failure signature maps $L_{\varphi} \colon \mathbb{R}^{D_\O} \to \mathbb{R}^{n}$ are:
\begin{eqnarray}\label{eqLOiMultiOutput}
L_{\varphi} &=& \left[
          \begin{array}{c}
            -\delta_{\varphi,1}\\
            \vdots\\
            -\delta_{\varphi,n_{S}}\\
            \textbf{0}_{(n - n_{S}) \times D_\O}\\
          \end{array}
        \right],
\end{eqnarray}
where $\delta_{\varphi,i} \in (\mathbb{R}_0^+)^{D_\O}$ and each component $\delta_{\varphi,i}(d)$ is the perturbation introduced by a failure ${\varphi}$ in the paths of $G_\O$ terminating with node $v_i$ and characterized by delay $d$. The following theorem states necessary and sufficient conditions to solve the EFPRG when $G_\O$ is multi-hop and $G_\R$ is single-hop.

\begin{theorem}\label{thMCNGRsingleGOmulti}
Let a MCN $N$ and the corresponding faulty set $\Phi$ be given, where $G_\R$ is single-hop and $G_\O$ is multi-hop with $n_{S} \geq n_{\Phi}$ terminating nodes. Then the EFPRG can be solved for each $\varphi \in \Phi$ if and only if the following condition holds:
\begin{equation}\label{conditionMCNGRsingleGOmulti}
d(\L_{\Phi}) = n_{\Phi},
\end{equation}
where the matrix ${L}_{\Phi} := \left[\begin{array}{cccc}
                     {L}_{\varphi_1} & {L}_{\varphi_2} & \cdots & {L}_{\varphi_{|\Phi|}}
                   \end{array}\right]$ is the juxtaposition of all failure signature maps in $\Phi$ and has dimensions $n_S \times n_\Phi$.

\begin{proof}
We need to state the equivalence between \eqref{conditionMCNGRsingleGOmulti} and \eqref{conditionExtFPRG}. For any $\varphi \in \Phi$, ${\L}_{\varphi} \subseteq \big(\mathcal N(C)\big)^\perp$ and $\bar{\L}_{\varphi} \subseteq \big(\mathcal N(C)\big)^\perp$. Thus, Lemma \ref{lemPropUOSA} implies that:
\begin{equation*}
\S^*\left(\bar{\L}_{\varphi}\right) = \bar{\L}_{\varphi} + \mathcal K_\varphi, \mbox{ } \mathcal K_\varphi \subseteq \mathcal N(C).
\end{equation*}
Moreover, for any $\varphi \in \Phi$, $\L_{\varphi} \cap \mathcal K_\varphi = \textbf{\emph{0}}$, thus:
\begin{equation*}
\S^*\left(\bar{\L}_{\varphi}\right) \cap \L_{\varphi} = \left(\bar{\L}_{\varphi} + \mathcal K_\varphi \right) \cap \L_{\varphi} = \bar{\L}_{\varphi} \cap \L_{\varphi}.
\end{equation*}
It follows that \eqref{conditionExtFPRG} is equivalent to the following:
\begin{equation}\label{conditionExtFPRGEquivalent}
\bar{\L}_{\varphi} \cap \L_\varphi = \textbf{\emph{0}}.
\end{equation}
Since $n_{S} \geq n_{\Phi}$ by assumption, then $d(\L_{\Phi}) \leq n_{\Phi}$. Since $\L_\varphi$ are monic, Condition \eqref{conditionExtFPRGEquivalent} implies that \eqref{conditionExtFPRG} holds if and only if $d(\L_{\Phi}) = n_{\Phi}$.
\end{proof}
\end{theorem}
\smallskip
The following theorem characterizes the relation between Condition \eqref{conditionMCNGRsingleGOmulti} and the topology of $G_\O(\eta_{\O})$.
\begin{theorem}\label{thGraphisTree}
Let a MCN $N$ and the corresponding faulty set $\Phi$ be given, where $G_\R$ is single-hop and $G_\O$ is multi-hop with $n_{S}$ terminating nodes. Then, $d(\L_{\Phi}) = n_{\Phi}$ if and only if $G_\O(\eta_{\O})$ is a tree, where $v_y$ is the root node and $v_1, \ldots, v_{n_{S}}$ are the leaves.

\begin{proof}
(\textbf{sufficiency}) Let $G_\O(\eta_{\O})$ be a tree, where~$v_y$ is the root node and the terminating nodes $v_1, \ldots, v_{n_{S}}$ are the leaves. Therefore, for each terminating node $v_i, i \in \{1, \ldots, n_S\}$ there exist a unique a link $e_i = (v_i',v_i)~\in~E_\O$, with $v_i' \in V_\O \setminus \{v_1, \ldots, v_{n_{S}}\}$. Define the configurations of failures $f_i = \{e_i\}, i \in \{1,\ldots,n_{S}\}$ and the corresponding failure signature maps $\{L_{f_1}, \ldots, L_{f_{n_S}}\}$, each characterized by $n_{S}$ rows and 1 column. Since $G_\O(\eta_{\O})$ is a tree, for each set $f \in 2^{E_\O} \setminus \big\{f_1, \ldots, f_{n_S}\big\}$, there exist $p \leq n_S$ and $e_1, \ldots, e_p$ such that $L_{f} m_{f}(kT) = \sum_{i=1}^p  L_{f_i} m_{f_i}(kT), \forall k \geq 0$. Since $\L_{f_i} \cap \L_{f_j} \neq \textbf{\emph{0}}$ for each $i,j = 1,\ldots,n_S$, $i \neq j$, then $\Phi = \{f_1, \ldots, f_{n_S}\}$ and $n_{\Phi} = n_{S}$. Since $L_{f_1}, \ldots, L_{f_{n_S}}$ are monic, then $d(\L_{\Phi}) = n_{\Phi}$.

(\textbf{necessity})
Assume that $G_\O(\eta_{\O})$ is not a tree. Then there exist nodes $v, \ v'$, and $v''$ such that $e'~=~(v',v),\ e''~=~(v'',v) \in E_\O$. Define $f' = \{e'\}$ and $f''=\{e''\}$ In this case, $L_{f'}$ assumes the following form:
\begin{equation*}
L_{f'} = -\left[
  \begin{array}{ccc}
    \delta'_{v_y,v_1}(D_\O) & \cdots & \delta'_{v_y,v_1}(1) \\
    \vdots & \ddots & \vdots \\
    \delta'_{v_y,v_{n_S}}(D_\O) & \cdots & \delta'_{v_y,v_{n_S}}(1) \\
    \textbf{0}_{(n - n_{S}) \times 1} & \cdots & \textbf{0}_{(n - n_{S}) \times 1}
  \end{array}
\right],
\end{equation*}
where $\delta'_{v_y,v_i}(d)$ is the contribution on the dynamics \eqref{eqFaultyMCNDynamicsSimultaneousFailures} of all paths starting from $v_y$, terminating in node $v_i$, passing through $e'$, and characterized by a delay $d$. It follows that:
\begin{equation*}
\L_{f'} \supseteq span\left(\left[
  \begin{array}{c}
    \sum_{d=1}^{D_\O} \delta'_{v_y,v_1}(d)\\
    \vdots\\
    \sum_{d=1}^{D_\O} \delta'_{v_y,v_{n_S}}(d)\\
    \textbf{0}_{(n - n_{S}) \times 1}
  \end{array}
\right]\right)
\end{equation*}
If a failure occurs in link $e'$, then the contribution $\sum_{d=1}^{D_\O} \delta'_{v_y,v_i}(d)$ on the dynamics \eqref{eqFaultyMCNDynamicsSimultaneousFailures} can be decomposed as the product of the contributions of all paths starting in $v_y$ and terminating in $v$ passing through $e'$, and of the contributions of all paths starting in $v$ and terminating in $v_i$. Thus,
\begin{equation*}
\L_{f'} \supseteq span\left(\left[
  \begin{array}{c}
    \left(\sum_{d=1}^{D_\O} \delta'_{v_y,v}(d)\right) \left(\sum_{d=1}^{D_\O} \delta_{v,v_1}(d)\right)\\
    \vdots\\
    \left(\sum_{d=1}^{D_\O} \delta'_{v_y,v}(d)\right) \left(\sum_{d=1}^{D_\O} \delta_{v,v_{n_S}}(d)\right)\\
    \textbf{0}_{(n - n_{S}) \times 1}
  \end{array}
\right]\right).
\end{equation*}
Since $L_{f''}$ can be defined similarly, then:
\begin{equation*}
\L_{f''} \supseteq span\left(\left[
  \begin{array}{c}
    \left(\sum_{d=1}^{D_\O} \delta''_{v_y,v}(d)\right) \left(\sum_{d=1}^{D_\O} \delta_{v,v_1}(d)\right)\\
    \vdots\\
    \left(\sum_{d=1}^{D_\O} \delta''_{v_y,v}(d)\right) \left(\sum_{d=1}^{D_\O} \delta_{v,v_{n_S}}(d)\right)\\
    \textbf{0}_{(n - n_{S}) \times 1}
  \end{array}
\right]\right).
\end{equation*}
It is clear that $\L_{f'} \cap \L_{f''} \neq \textbf{\emph{0}}$. If $\exists k \geq 0 : L_{f'} m_{f'}(kT) \neq L_{f''} m_{f''}(kT)$, then the configurations of failures $f'$ and $f''$ belong to different equivalence classes of $\Phi$ and thus $d(\L_{\Phi}) < n_{\Phi}$. If $L_{f'} m_{f'}(kT) = L_{f''} m_{f''}(kT), \forall k \geq 0$, then the configurations of failures of $f'$ and $f''$ belong to the same equivalence class $[L_{f'} m_{f'}]$ of $\Phi$, and we can not conclude that $d(\L_{\Phi}) < n_{\Phi}$. However, the simultaneous failure of links $e'$ and $e''$ belongs to the equivalence class $[L_{f' \cup f''} m_{f' \cup f''}]$, with $L_{f' \cup f''} \neq L_{f'}$ and $\L_{f' \cup f''}~\cap~\L_{f'} \neq \textbf{\emph{0}}$, and thus $d(\L_{\Phi}) < n_{\Phi}$.
\end{proof}
\end{theorem}
\begin{corollary}\label{corLOisKerCperp}
Let a MCN $N$ and the corresponding faulty set $\Phi$ be given, where $G_\R$ is single-hop and $G_\O$ is multi-hop with $n_{S}$ terminating nodes. If the EFPRG can be solved for each $\varphi \in \Phi$, then $n_{S} = n_{\Phi}$ and $\L_{\Phi} = \big(\mathcal N(C)\big)^\perp$.

\begin{proof}
Straightforward since $G_\O(\eta_{\O})$ is a tree, and thus to each terminating node $v_i, i \in \{1, \ldots, n_S\}$ corresponds only one path from $v_y$ to $v_i$.
\end{proof}
\end{corollary}
The necessary and sufficient condition given in Theorem \ref{thGraphisTree} provides a hard constraint on the topology of $G_\O(\eta_{\O})$ induced by the scheduling $\eta_{\O}$. This is not surprising, since we require to solve the EFPRG for the set $\Phi$ of \emph{all} configurations of failures that perturb the dynamics \eqref{eqFaultyMCNDynamicsSimultaneousFailures}. From an implementation point of view, this constraint can be both interpreted as hardware or software redundancy. In the former case, the tree structure of $G_\O(\eta_{\O})$ provides a hardware separation for all paths from $v_y$ to the terminating nodes. However, a tree communication graph might be not always implementable in real cases: therefore, the constraint on $G_\O(\eta_{\O})$ can be implemented by using, for those communication nodes that receive data from multiple incoming links, separate memory slots for each of the incoming data. These nodes will transmit distinct data for each memory slot, thus providing a software separation for all paths from $v_y$ to the terminating nodes. In general, a combination of the above approaches is reasonably implementable in a real communication network. An interesting future research direction is relating the properties of $G_\O(\eta_{\O})$ with Condition \eqref{conditionMCNGRsingleGOmulti} when the number of simultaneous failures that can occur is bounded, or when failures can not occur in some \emph{secure} paths of the communication network.


\subsection{$G_\R$ and $G_\O$ multi-hop} \label{subsecGRGOmultihop}

\begin{figure*}[t]
\begin{center}
\includegraphics[width=30px]{blank}
\end{center}
\end{figure*}

\begin{figure*}[t]
\begin{picture}(0,0)
\put(0,30){
\scalebox{0.65}{
$
\left[
  \begin{array}{c|cc|ccccccc}
    \textbf{0}_{1 \times n_{S}-1} & 0 & 0         & \gamma_1(1) C_\P B_\P         & \sum\limits_{i=1}^2 \gamma_1(i) C_\P A_\P^{2-i} B_\P              & \cdots & \sum\limits_{i=1}^{D_\O-1} \gamma_1(i) C_\P A_\P^{(D_\O-1-i)} B_\P                       & \sum\limits_{i=1}^{D_\O} \gamma_1(i) C_\P A_\P^{(D_\O-i)} B_\P                            & \sum\limits_{i=1}^{D_\O} \gamma_1(i) C_\P A_\P^{(D_\O+1-i)} B_\P                            & \cdots\\
    \hline
    & 0 & 0         & \gamma_{2}(1) C_\P B_\P  & \sum\limits_{i=1}^2 \gamma_{2}(i) C_\P A_\P^{2-i} B_\P       & \cdots & \sum\limits_{i=1}^{D_\O-1} \gamma_{2}(i) C_\P A_\P^{(D_\O-1-i)} B_\P                & \sum\limits_{i=1}^{D_\O} \gamma_{2}(i) C_\P A_\P^{(D_\O-i)} B_\P                     & \sum\limits_{i=1}^{D_\O} \gamma_{2}(i) C_\P A_\P^{(D_\O+1-i)} B_\P                     & \cdots\\
    \textbf{I}_{n_s-1} & \vdots & \vdots    & \vdots                        & \vdots                                                            & \ddots & \vdots                                                                                   & \vdots                                                                                    & \vdots                                                                                        & \cdots\\    & 0 & 0         & \gamma_{n_{S}}(1) C_\P B_\P  & \sum\limits_{i=1}^2 \gamma_{n_{S}}(i) C_\P A_\P^{2-i} B_\P       & \cdots & \sum\limits_{i=1}^{D_\O-1} \gamma_{n_{S}}(i) C_\P A_\P^{(D_\O-1-i)} B_\P                & \sum\limits_{i=1}^{D_\O} \gamma_{n_{S}}(i) C_\P A_\P^{(D_\O-i)} B_\P                     & \sum\limits_{i=1}^{D_\O} \gamma_{n_{S}}(i) C_\P A_\P^{(D_\O+1-i)} B_\P                     & \cdots\\
    \hline
    \textbf{0}_{1 \times n_{S}-1} & 0 & 0         & 0                             & 0                                                                 & \cdots & C_\P B_\P                                                                                & C_\P A_\P B_\P                                                                            & C_\P A_\P^2 B_\P                                                                              & \cdots\\
    \textbf{0}_{1 \times n_{S}-1} & 0 & 0         & 0                             & 0                                                                 & \cdots & C_\P A_\P B_\P                                                                           & C_\P A_\P^2 B_\P                                                                          & C_\P A_\P^3 B_\P                                                                              & \cdots\\
    \vdots & \vdots & \vdots    & \vdots                        & \vdots                                                            & \ddots & \vdots                                                                                   & \vdots                                                                                    & \vdots                                                                                        & \ddots\\
    \textbf{0}_{1 \times n_{S}-1} & 0 & 0         & 0                             & 0                                                                 & \cdots & C_\P A_\P^{D_\O-4} B_\P                                                                  & C_\P A_\P^{D_\O-3} B_\P                                                                        & C_\P A_\P^{D_\O-2} B_\P                                                                            & \cdots\\
    \textbf{0}_{1 \times n_{S}-1} & 0 & 0         & 0                             & C_\P B_\P                                                         & \cdots & C_\P A_\P^{D_\O-3} B_\P                                                                  & C_\P A_\P^{D_\O-2} B_\P                                                                       & C_\P A_\P^{D_\O-1} B_\P                                                                            & \cdots\\
    \textbf{0}_{1 \times n_{S}-1} & 0 & 0         & C_\P B_\P                     & C_\P A_\P B_\P                                                    & \cdots & C_\P A_\P^{D_\O-2} B_\P                                                                  & C_\P A_\P^{D_\O-1} B_\P                                                                       & C_\P A_\P^{D_\O} B_\P                                                                              & \cdots\\
    \textbf{0}_{1 \times n_{S}-1} & 0 & B_\P      & A_P B_P                       & A^2_P B_P                                                         & \cdots & A_P^{D_\O-1} B_P                                                                         & A^{D_\O}_P B_P                                                                            & A^{D_\O+1}_P B_P& \cdots\\
    \textbf{0}_{1 \times n_{S}-1} & 1 & 0 & 0 & 0 & \cdots & 0 & 0 & 0 & \cdots\\
    \textbf{0}_{n_{\R}-1 \times n_{S}-1} & \textbf{0}_{n_{\R}-1 \times 1} & \textbf{0}_{n_{\R}-1 \times 1} & \textbf{0}_{n_\R-1 \times 1} & \textbf{0}_{n_\R-1 \times 1} & \cdots & \textbf{0}_{n_\R-1 \times 1} & \textbf{0}_{n_\R-1 \times 1} & \textbf{0}_{n_\R-1 \times 1} & \cdots
  \end{array}
\right]
$
}
}
\end{picture}
\vspace{1.5cm}
\caption{Inductive definition of matrix $\Psi_\infty$.} \label{figMatrixH}
\end{figure*}

When both $G_\R$ and $G_\O$ are multi-hop, we need to define the set $\Phi = \Phi_{\R} \cup \Phi_{\O}$ of equivalence classes that equally perturb the dynamics \eqref{eqFaultyMCNDynamicsSimultaneousFailures}. In this case, failures occur in both the controllability and observability graphs. Therefore, by an appropriate choice of $m_\varphi(kT)$, we define the failure signature maps associated to the equivalence classes $\varphi_\R \in \Phi_\R$ and $\varphi_\O \in \Phi_\O$ by:
\begin{equation*}
L_{\varphi_\R} = \left[
          \begin{array}{l}
            \textbf{0}_{(n_\O + n_\P) \times n_\R} \\
            -\delta_{\varphi_\R}\\
            \textbf{0}_{(n_\R -1) \times n_\R} \\
          \end{array}
        \right], \mbox{ }
L_{\varphi_\O} = \left[
          \begin{array}{lll}
            -\delta_{\varphi_\O}\\
            \textbf{0}_{(n - n_S) \times n_\O}\\
          \end{array}
        \right],
\end{equation*}
with $\delta_{\varphi_\R} \in (\mathbb{R}_0^+)^{D_\R}$ a row vector, and $\delta_{\varphi_\O} \in (\mathbb{R}_0^+)^{n_S \times D_\O}$ as defined in \eqref{eqLOiMultiOutput}.

We recall that, for each $\varphi_\R \in \Phi_\R$ non-empty, $\L_{\varphi_\R}~=~span(\textbf{e}_{n_\O + n_\P + 1})$. Therefore, we will consider w.l.o.g. only one failure in the reachability graph, namely $\Phi_\R = \{\varnothing, \varphi_\R\}$ with $\L_{\varphi_\R}~=~span(\textbf{e}_{n_\O + n_\P + 1})$.

Moreover, by Theorem \ref{thGraphisTree}, a necessary condition to solve the EFPRG for any $\varphi_\O \in \Phi_\O$ is that $G_\O$ is a tree. Therefore, we will consider w.l.o.g. a failure in the observability graph for each path, namely $\Phi_\O = \{\varphi_{1}, \ldots, \varphi_{{n_S}}\}$ with $L_{\varphi_{i}} = span(\textbf{e}_{i})$.

The following theorem states that it is not possible to detect failures in the controllability and observability graphs using the measurements of the observability graph.
\begin{theorem}\label{thMCNGRmultiGOmulti}
Let a MCN $N$ and the corresponding faulty set $\Phi$ be given, where $G_\R$ is multi-hop and $G_\O$ is multi-hop with $n_{S}$ terminating nodes. Then the EFPRG is not solvable for any $\varphi_\R \in \Phi_\R$ and any $\varphi_\O \in \Phi_\O$.

\begin{proof}
We first show that $\S^*\left(\bar\L_{\varphi_\R}\right) \cap \L_{\varphi_\R} \neq \textbf{\emph{0}}$. By Corollary \ref{corLOisKerCperp}, $\sum_{\varphi_\O \in \Phi_\O}\L_{\varphi_\O} = \big(\mathcal N(C)\big)^\perp$, and $\S^*(\sum_{\varphi_\O \in \Phi_\O}\L_{\varphi_\O})~=~\mathbb{R}^n$ by Lemma \ref{lemPropUOSA}. Since $\bar\L_{\varphi_\R} = \sum_{\varphi_\O \in \Phi_\O}\L_{\varphi_\O}$, then $\S^*\left(\bar\L_{\varphi_\R}\right)~\cap~ \L_{\varphi_\R} \neq \textbf{\emph{0}}$.

%

To complete the proof, we need to show that for each $i \in \{1, \ldots, n_S\}$, $\S^*\left(\bar{\L}_{\varphi_{{i}}}\right) \cap \L_{\varphi_{{i}}} \neq \textbf{\emph{0}}$, with $\varphi_{{i}} \in \Phi_\O$. We will only provide the proof for $i=1$: the same reasoning can be used for $i \in \{2, \ldots, n_S\}$.

The space $\W^*\left(\bar{\L}_{\varphi_{{1}}}\right)$ is generated by the submatrix $\Psi_{h}$, which consists of the first $h$ columns of the matrix $\Psi_\infty$ with infinite columns inductively defined in Figure \ref{figMatrixH}, and where the value of $h$ depends on the terminating condition of the CAISA Algorithm. More precisely, $h$ is the smallest integer such that $rank\big(span(\Psi_{h}) \cap \mathcal N(C)\big) = rank\big(span(\Psi_{h+1}) \cap \mathcal N(C)\big)$. The above terminating condition occurs at column $h$ if and only if one of the following two conditions holds: (i) the $1$-st row of column $h$ (which is a scalar) is equal to zero and column $h$ is linearly dependent on all the previous columns $1, \ldots, h-1$; (ii) the $1$-st row of column $h$ is different from zero. We show in the following that condition (ii) will always stop the CAISA algorithm before condition (i) can occur.

Let $m \in \mathbb N \cup \{0\}$ be the smallest value such that $C_\P A_\P^m B_\P \neq 0$. Since $(A_\P,B_\P)$ is controllable and $(C_\P, A_\P)$ is observable, then $m \leq n_\P-1$. Note that the first $n_s + 1$ columns of $\Psi_\infty$ are already present, since they belong to $\bar{\L}_{\varphi_{{1}}}$. The subsequent $m$ columns are linearly independent from the previous columns since $(A_\P,B_\P)$ is controllable and $m \leq n_\P-1$. Since the scalar $C_\P A_\P^m B_\P \neq 0$ appears at row $n_S + D_\O - 1$ and at column $n_S + 2 + m$, the subsequent $D_\O-2$ columns are linearly independent from the previous columns. Therefore, column $h$ can be linearly dependent on all the previous columns for $h \geq h_1 = n_S + m + D_\O +1$.

Let $1 \leq d_{1} \leq D_\O$ be the smallest value such that $\gamma_{1}(d_{1}) \neq 0$. Therefore, the $1$-st row of $\Psi_\infty$ will have a non-zero value for the first time at row column $h_2 = n_S + m + d_{1} +1$. Since $h_2 \leq h_1$, then condition (ii) will always stop the CAISA algorithm before condition (i) can occur. Therefore:
$$
\W^*\left(\bar{\L}_{\varphi_{{1}}}\right) = span\left[\begin{array}{ccc}
                                                      \textbf{I}_{n_S-1} & \textbf{0} & \psi_{1} \\
                                                      0 & 0 & \psi_{2} \\
                                                      \textbf{0} & \textbf{I}_{l} & \psi_{3}
                                                    \end{array}\right],
$$
where $l \leq n - n_S$, $\psi_{1}$ is a $n_S -1$ column vector, $\psi_{2} \neq 0$ is a scalar, and $\psi_{3}$ is a $l$ column vector. Applying the UOSA algorithm, we obtain:
$$
\S_1\left(\bar{\L}_{\varphi_{{1}}}\right) = \W^*\left(\bar{\L}_{\varphi_{{1}}}\right) + \mathcal N(C) = \mathbb R^n = \S^*\left(\bar{\L}_{\varphi_{{1}}}\right),
$$
which clearly implies that $\S^*\left(\bar{\L}_{\varphi_{{1}}}\right) \cap \L_{\varphi_{{1}}} \neq 0$.
\end{proof}
\end{theorem}

Theorem \ref{thMCNGRmultiGOmulti} states that, in order to detect failures in the observability graph, the controllability graph must not be subject to failures. By a practical point of view, the communication protocol in the controllability graph is required to implement failure detection using handshaking messages between nodes and inform the controller about the set of faulty links.

\bibliographystyle{IEEEtran}
\bibliography{mcnbib}

\begin{thebibliography}{10}
\providecommand{\url}[1]{#1}
\csname url@rmstyle\endcsname
\providecommand{\newblock}{\relax}
\providecommand{\bibinfo}[2]{#2}
\providecommand\BIBentrySTDinterwordspacing{\spaceskip=0pt\relax}
\providecommand\BIBentryALTinterwordstretchfactor{4}
\providecommand\BIBentryALTinterwordspacing{\spaceskip=\fontdimen2\font plus
\BIBentryALTinterwordstretchfactor\fontdimen3\font minus
  \fontdimen4\font\relax}
\providecommand\BIBforeignlanguage[2]{{%
\expandafter\ifx\csname l@#1\endcsname\relax
\typeout{** WARNING: IEEEtran.bst: No hyphenation pattern has been}%
\typeout{** loaded for the language `#1'. Using the pattern for}%
\typeout{** the default language instead.}%
\else
\language=\csname l@#1\endcsname
\fi
#2}}

\bibitem{akyildiz_wireless_2004}
{I.F. Akyildiz} and {I.H. Kasimoglu}, ``{Wireless Sensor and Actor Networks:
  Research Challenges},'' \emph{{Ad Hoc Networks}}, vol.~2, no.~4, pp.
  351--367, 2004.

\bibitem{song_wirelesshart:_2008}
J.~Song, S.~Han, {A.K. Mok}, D.~Chen, M.~Lucas, M.~Nixon, and W.~Pratt,
  ``{WirelessHART: Applying Wireless Technology in Real-Time Industrial Process
  Control},'' in \emph{{RTAS}}, 2008.

\bibitem{song_complete_2008}
J.~Song, S.~Han, X.~Zhu, {A.K. Mok}, D.~Chen, and M.~Nixon, ``{A Complete
  WirelessHART Network},'' in \emph{ACME}, 2008, pp. 381--382.

\bibitem{Zhang2001}
{W. Zhang}, {M.S. Branicky}, and {S.M. Phillips}, ``{Stability of Networked
  Control Systems},'' \emph{IEEE Control Systems Magazine}, vol.~21, no.~1, pp.
  84--99, February 2001.

\bibitem{WalshCSM2001}
{G.C. Walsh} and H.~Ye, ``{Scheduling of Networked Control Systems},''
  \emph{IEEE Control Systems Magazine}, pp. 57--65, February 2001.

\bibitem{SpecialIssueNCS2004}
P.~Antsaklis and J.~Baillieul, ``{Guest Editorial Special Issue on Networked
  Control Systems},'' \emph{IEEE Transactions on Automatic Control}, vol.~49,
  no.~9, pp. 1421--1423, September 2004.

\bibitem{Hespanha2007}
{J.P. Hespanha}, P.~Naghshtabrizi, and Y.~Xu, ``{A Survey of Recent Results in
  Networked Control Systems},'' \emph{Proceedings of the IEEE}, vol.~95, no.~1,
  pp. 138--162, January 2007.

\bibitem{HeemelsTAC10}
{W.P.M.H. Heemels}, {A.R. Teel}, {N. van de Wouw}, and {D. Ne\v{s}i\'{c}},
  ``{Networked Control Systems With Communication Constraints: Tradeoffs
  Between Transmission Intervals, Delays and Performance},'' \emph{IEEE
  Transactions on Automatic Control}, vol.~55, no.~8, pp. 1781 --1796, August
  2010.

\bibitem{andersson_simulation_2005}
M.~Andersson, D.~Henriksson, A.~Cervin, and K.~Arzen, ``{Simulation of Wireless
  Networked Control Systems},'' in \emph{{Proceedings of the 44th IEEE
  Conference on Decision and Control and European Control Conference}}, 2005,
  pp. 476--481.

\bibitem{AlurRTAS09}
R.~Alur, {A. D'Innocenzo}, {K.H. Johansson}, {G.J. Pappas}, and {G. Weiss},
  ``{Modeling and Analysis of Multi-Hop Control Networks},'' in
  \emph{Proceedings of the 15th IEEE Real-Time and Embedded Technology and
  Applications Symposium (RTAS)}, 2009.

\bibitem{AlurTAC11}
{R. Alur}, {A. D'Innocenzo}, {K.H. Johansson}, {G.J. Pappas}, and {G. Weiss},
  ``{Compositional Modeling and Analysis of Multi-Hop Control Networks},''
  \emph{IEEE Transactions on Automatic Control}, 2011, accepted for publication
  as regular paper.

\bibitem{DiBenedettoIFAC11Stab}
{M.D. Di Benedetto}, A.~D'Innocenzo, and E.~Serra, ``{Fault Tolerant
  Stabilizability of Multi-Hop Control Networks},'' in \emph{Proceedings of the
  18th IFAC World Congress, Milan, Italy}, 2011, preprint available at
  arXiv:1103.4340v1.

\bibitem{D'InnocenzoCASE2009}
A.~D'Innocenzo, G.~Weiss, R.~Alur, {A.J. Isaksson}, {K.H. Johansson}, and {G.J.
  Pappas}, ``{Scalable Scheduling Algorithms for Wireless Networked Control
  Systems},'' in \emph{Proceedings of the 5th IEEE Conference on Automation
  Science and Engineering (CASE)}, 2009.

\bibitem{ZavlanosCDC07}
{M.M. Zavlanos} and {G.J. Pappas}, ``{Distributed Connectivity Control of
  Mobile Networks},'' in \emph{{Proceedings of the 46th IEEE Conference on
  Decision and Control}}, December 2007, pp. 3591 --3596.

\bibitem{MeskinTAC2009}
N.~Meskin and K.~Khorasani, ``{Actuator Fault Detection and Isolation for a
  Network of Unmanned Vehicles},'' \emph{{IEEE Transactions on Automatic
  Control}}, vol.~54, no.~4, pp. 835 --840, April 2009.

\bibitem{BeardPhDThesis1971}
R.~Beard, ``{Failure Accomondation in Linear Systems Through
  Self-Reorganization},'' Ph.D. dissertation, {MIT}, 1971.

\bibitem{JonesPhDThesis1973}
H.~Jones, ``{Failure Detection in Linear Systems},'' Ph.D. dissertation, {MIT},
  1973.

\bibitem{MassoumniaTAC89}
{M.-A. Massoumnia}, {G.C. Verghese}, and {A.S. Willsky}, ``{Failure Detection
  and Identification},'' \emph{IEEE Transactions on Automatic Control},
  vol.~34, no.~3, pp. 316 --321, Mar. 1989.

\bibitem{DePersisTAC2001}
{C. De Persis} and A.~Isidori, ``{A Geometric Approach to Nonlinear Fault
  Detection and Isolation},'' \emph{{IEEE Transactions on Automatic Control}},
  vol.~46, no.~6, pp. 853 --865, June 2001.

\bibitem{Gupta10}
R.~Gupta and M.-Y. Chow, ``{Networked Control System: Overview and Research
  Trends},'' \emph{IEEE Transactions on Industrial Electronics}, vol.~57,
  no.~7, pp. 2527 --2535, July 2010.

\bibitem{WangTSP2008}
Y.~Wang, {S.X. Ding}, H.~Ye, and G.~Wang, ``{A New Fault Detection Scheme for
  Networked Control Systems Subject to Uncertain Time-Varying Delay},''
  \emph{{IEEE Transactions on Signal Processing}}, vol.~56, no.~10, pp. 5258
  --5268, October 2008.

\bibitem{CommaultTAC2007}
C.~Commault and {J.-M. Dion}, ``{Sensor Location for Diagnosis in Linear
  Systems: A Structural Analysis},'' \emph{{IEEE Transactions on Automatic
  Control}}, vol.~52, no.~2, pp. 155 --169, February 2007.

\bibitem{PappasCDC2010a}
{S. Sundaram}, {M. Pajic}, {C.N. Hadjicostis}, {R. Mangharam}, and {G.J.
  Pappas}, ``{The Wireless Control Network: Monitoring for Malicious
  Behavior},'' in \emph{Proceedings of the 49th IEEE Conference on Decision and
  Control (CDC)}, December 2010, pp. 5979 --5984.

\bibitem{DiBenedettoIFAC11Rout}
{M.D. Di Benedetto}, A.~D'Innocenzo, and E.~Serra, ``{Dynamical Power
  Optimization by Decentralized Routing Control in Multi-Hop Wireless Control
  Networks},'' in \emph{Proceedings of the 18th IFAC World Congress, Milan,
  Italy}, 2011.

\bibitem{WonhamLMC1979}
{W.M. Wonham}, \emph{{Linear Multivariable Control: a Geometric Approach}},
  2nd~ed., ser. {Applications of Mathematics}.\hskip 1em plus 0.5em minus
  0.4em\relax Springer-Verlag, 1979.

\end{thebibliography}
\end{document}